\documentclass{amsart}
\usepackage{amsmath,amscd,amssymb,amsthm}

\newtheorem{lemma}{Lemma}[section]

\newtheorem{prop}[lemma]{Proposition}
\newtheorem{thm}[lemma]{Theorem}

\newtheorem{conj}{Conjecture}

\def\Uhelt{h}
\def\HC{\mathop{\rm HC}\nolimits}
\def\rank{\mathop{\rm rank}\nolimits}

\def\cal{\mathcal}

\def\I{{\cal I}}

\def\O{{\cal O}}

\def\frak{\mathfrak}

\def\da{{\frak a}}

\def\dF{{\frak F}}

\def\dg{{\frak g}}
\def\dgl{\dg\dl}
\def\dh{{\frak h}}

\def\dk{{\frak k}}
\def\dl{{\frak l}}
\def\dm{{\frak m}}
\def\dn{{\frak n}}
\def\dgo{{\frak o}}

\def\ds{{\frak s}}
\def\dsl{\ds\dl}

\def\du{{\frak u}}
\def\dU{{\frak U}}

\def\dz{{\frak z}}

\def\dZ{{\frak Z}}

\def\Bbb{\mathbb}

\def\bC{\Bbb C}

\def\bN{\Bbb N}

\def\bR{\Bbb R}

\def\bZ{\Bbb Z}

\def\ep{\epsilon}

\def\la{\langle}

\def\ra{\rangle}

\def\t{\tilde}

\def\End{\mathop{\rm End}\nolimits}
\def\Frac{\mathop{\rm Frac}\nolimits}

\def\oh{{\ts\frac{1}{2}}}
\def\scrm{\scriptsize\rm}

\def\Span{\mathop{\rm Span}\nolimits}
\def\ss{{\mbox{\scrm ss}}}

\def\thup{{\mbox{\scrm th}}}
\def\ts{\textstyle}

\def\bgno{\bigbreak\noindent}

\def\expro{extremal projector}

\def\hwv{highest weight vector}

\def\ie{{\em i.e.,\/}}
\def\iff{if and only if}
\def\irr{irreducible}

\def\meno{\medbreak\noindent}

\def\r{representation}

\def\uea{universal enveloping algebra}

\title[Factorizations of relative extremal projectors]{Factorizations of relative extremal projectors}

\author{Charles H.\ Conley}
\address{Department of Mathematics \\University of North Texas \\Denton TX 76203, USA} 
\email{conley@unt.edu}

\author{Mark R.\ Sepanski}
\address{Department of Mathematics \\Baylor University \\Waco TX 76798, USA}
\email{mark\_sepanski@baylor.edu}

\thanks{The first author was partially supported by Simons Foundation Collaboration Grant 207736.}

\begin{document}

\begin{abstract}
We survey earlier results on factorizations of extremal projectors and relative extremal projectors and present preliminary results on non-commutative factorizations of relative extremal projectors: we deduce the existence of such factorizations for $\dsl_4$ and $\dsl_5$.

\bgno {\sc 2010 Mathematics Subject Classifications:} 17B20, 17B35

\bgno {\em This article is dedicated to V.\ S.\ Varadarajan on the occasion of his retirement.\/}

\end{abstract}

\maketitle

%%%%%%%%%%%%%%%%
\section{Introduction}  \label{Intro} 
%%%%%%%%%%%%%%%%

Extremal projectors were first investigated in the 1960's and 1970's by Asherova, Smirnov, and Tolstoi; their results are summarized in \cite{AST79}.  Zhelobenko wrote a series of articles on projectors in the 1980's and 1990's, including the survey monograph \cite{Zh90}.  In this note we will focus on relative extremal projectors, but we begin with a brief account of extremal projectors.  We will take the liberty of attributing to \cite{AST79} and \cite{Zh90} some results which were in fact first discovered in earlier works of the same authors: see the references of those papers.

Given any Lie algebra $\dk$, we have the \uea\ $\dU(\dk)$ and the augmentation ideal $\dU^+(\dk) := \dk \dU(\dk)$.  If $W$ is a $\dk$-module, we write $W^\dk$ for the space of $\dk$-invariants in $W$.  Throughout this article we write $\bN$ for the non-negative integers and $\bZ^+$ for the positive integers.

Let $\dg$ be a finite dimensional complex reductive Lie algebra, $\dh$ a Cartan subalgebra, and $\dn^- \oplus \dh \oplus \dn^+$ a triangular decomposition of $\dg$.  We denote the associated positive and simple root systems by $\Delta(\dn^+)$ and $\Pi(\dn^+)$, respectively.  More generally, given any $\dh$-module $V$ we write $\Delta(V)$ for its weights in $\dh^*$ and $V_\mu$ for its $\mu$-weight space.  Throughout the article, all $\dh$-invariant subalgebras of $\dg$ will be endowed with the positive and negative systems inherited from $\dg$.  We sometimes write $\dg_\ss$ for the semisimple part of $\dg$ and $\dz(\dg)$ for its center, so that $\dg = \dg_\ss \oplus \dz(\dg)$.

One may ask the following na\"ive question: is there an element of $\dU(\dg)$ which projects any \r\ $V$ of $\dg$ in the category $\O(\dg)$ to its highest weight space $V^{\dn^+}$ along the sum of its lower weight spaces $\dn^- V$?  The answer is no, but there is such an element in a certain extension $\dF(\dg)$ of $\dU(\dg)$, the {\em extremal projector\/} $P(\dg)$.  Its action is defined on all weight spaces $V_\mu$ such that
\begin{equation} \label{Pg domain}
   V_\mu^{\dn^+} \cap \dn^- V = 0.
\end{equation}

In order to define $\dF(\dg)$, fix an $\dsl_2$-triple $\{E_\alpha, F_\alpha, H_\alpha\}$ for each positive root $\alpha$ in $\Delta(\dn^+)$.  Thus $E_\alpha$ spans $\dn^+_\alpha$, $F_\alpha$ spans $\dn^-_{-\alpha}$, $H_\alpha$ is the element $[E_\alpha, F_\alpha]$ of $\dh$, and $\alpha(H_\alpha) = 2$.  We index the positive roots and use multinomial notation:
\begin{equation} \label{roots}
   \Delta(\dn^+) := \bigl\{ \alpha_1, \ldots, \alpha_m \bigr\}, \quad
   E^I := E_{\alpha_1}^{I_1} \cdots E_{\alpha_m}^{I_m}, \quad
   F^I := F_{\alpha_1}^{I_1} \cdots F_{\alpha_m}^{I_m},
\end{equation}
where $I \in \bN^m$.  For any $K \in \bZ^m$, let $|K|$ denote the weight $\sum_{r=1}^m K_r \alpha_r$.

\meno {\bf Definition.} {\em
$\dF(\dg)$ is the algebraic direct sum\/ $\bigoplus_{\gamma \in \Delta(\dU(\dg))} \dF(\dg)_\gamma$ of its weight spaces, where\/ $\dF(\dg)_\gamma$ is the space of formal series in the monomials\/ $F^I E^J$ of weight\/ $\gamma$ with coefficients in the fraction field\/ $\Frac \dU(\dh)$ of\/ $\dU(\dh)$:\/}
\begin{equation*}
   \dF(\dg)_\gamma := \Bigl\{ \sum_{|J - I| = \gamma}
   F^I E^J \Uhelt_{IJ}: \Uhelt_{IJ} \in \Frac \dU(\dh) \Bigr\}.
\end{equation*}

The commutation relations of $\dU(\dg)$ extend to an algebra structure on $\dF(\dg)$.  Recall that the Cartan involution $\theta$ of $\dg$ exchanges $E_\alpha$ and $-F_\alpha$ for $\alpha$ simple and is $-1$ on $\dh$.  Let $\Omega \mapsto \Omega^*$ be the {\em Hermitian anti-involution\/} of $\dF(\dg)$, which is $-\theta$ on $\dg$.  Elements of $\dF(\dg)$ fixed by this anti-involution are called {\em Hermitian.\/}  For reference, note that $\dg = \dh$ gives
\begin{equation*}
   \dF(\dh) = \Frac \dU(\dh).
\end{equation*}

\begin{thm} \label{expro dfn} \cite{AST79}
There is a unique non-zero Hermitian idempotent $P(\dg)$ in $\dF(\dg)_0$, the\/ {\em \expro,} such that
\begin{equation*}
   \dn^+ P(\dg) = 0 = P(\dg) \dn^-.
\end{equation*}
\end{thm}

The {\it universal Verma module\/} $M(\dg)$ is $\dF(\dg) / \dF(\dg) \dn^+$.  It is a two-sided $\dF(\dh)$-module, spanned freely by the image of $\dU(\dn^-)$ under both the right and left actions.  Under the adjoint action of $\dh$, $F^I$ has weight $-|I|$ and $M(\dg)$ is the direct sum of its weight spaces:
\begin{equation} \label{h decomp of M}
   M(\dg) = \bigoplus_{\nu \in \Delta(\dU(\dn^+))} M(\dg)_{-\nu}, \qquad
   M(\dg)_{-\nu} = \Span_{\dF(\dh)} \bigl\{ F^I: |I| = \nu \bigr\}.
\end{equation}

Write $\End_{\rho(\dh)} M(\dg)$ for the endomorphisms of $M(\dg)$ commuting with the right action $\rho$ of $\dh$.  Since $\dF(\dg)$ acts on $M(\dg)$ from the left, there is a natural homomorphism from $\dF(\dg)$ to $\End_{\rho(\dh)} M(\dg)$.

The {\em Shapovalov form\/} $\la \cdot, \cdot \ra$ on $M(\dg)$ is $\dF(\dh)$-valued: for $\Theta$ and $\Theta'$ in the space $\dU(\dn^-) \dF(\dh)$, $\la \Theta, \Theta' \ra$ is the projection of the weight zero component of $\Theta^* \Theta'$ to $\dF(\dh)$ along $(\dn^- \dF(\dg) \dn^+)_0$.  It is non-degenerate, symmetric, and right $\dh$-bilinear, and Hermitian elements of $\dF(\dg)$ act on $M(\dg)$ by Hermitian operators.

The universal Verma module is the Hermitian direct sum of its highest weight space $M(\dg)^{\dn^+}$ and the sum of its lower weight spaces $\dn^- M(\dg)$:
\begin{align}
   & M(\dg)^{\dn^+} = M(\dg)_0, \nonumber \\[6pt]
   & \dn^- M(\dg) = \bigoplus_{\nu \in \Delta(\dU^+(\dn^+))} M(\dg)_{-\nu}, \label{Herm decomp of M} \\[6pt]
   & M(\dg) = M(\dg)^{\dn^+} \oplus \dn^- M(\dg). \nonumber
\end{align}

\begin{thm} \label{Verma} \cite{Zh90}
The homomorphism from $\dF(\dg)$ to\/ $\End_{\rho(\dh)} M(\dg)$ is an isomorphism.  It carries $P(\dg)$ to the projection from $M(\dg)$ to $M(\dg)^{\dn^+}$ along\/ $\dn^- M(\dg)$.
\end{thm}

Extremal projectors have applications to Mickelsson step algebras and Yangians; see \cite{Zh90} and the comprehensive text of Molev \cite{Mo07}.  They also arise in the work of Etingof, Tarasov, Varchenko, and others on dynamical quantum Weyl groups \cite{TV00, EV02}.

%%%%%%%%%%%%%%%%
\section{Factorizations of the extremal projector}  \label{Facs} 
%%%%%%%%%%%%%%%%

\subsection{Non-commutative factorizations}

As usual, write $w_0$ for the longest element of the Weyl group $W(\dg)$, $\rho_\dg$ for the half-sum $\oh \sum_{\Delta(\dn^+)} \alpha$ of the positive roots, and $s_\alpha$ for the reflection in a root $\alpha$.

\meno {\bf Definition.}
An ordering $(\alpha_1, \ldots, \alpha_m)$ of $\Delta(\dn^+)$ is {\em normal\/} if whenever $\alpha_r + \alpha_s$ is a root $\alpha_t$, then either $r < t < s$ or $s < t < r$.

\begin{prop} \cite{Zh87}
Normal orders are in bijection with reduced expressions for $w_0$: if $\alpha_1', \ldots, \alpha_m'$ are elements of\/ $\Pi(\dn^+)$ (usually not distinct) such that $w_0 = s_{\alpha_1'} \cdots s_{\alpha_m'}$, then $(\alpha_1, \ldots, \alpha_m)$ is a normal order of\/ $\Delta(\dn^+)$ for $\alpha_r = s_{\alpha_1'} \cdots s_{\alpha_{r-1}'} \alpha_r'$.

Conversely, if $(\alpha_1, \ldots, \alpha_m)$ is a normal order of\/ $\Delta(\dn^+)$, then the roots $\alpha_r' := s_{\alpha_1} \cdots s_{\alpha_{r-1}} \alpha_r$ are in $\Pi(\dn^+)$ and $w_0 = s_{\alpha_1'} \cdots s_{\alpha_m'}$. 
\end{prop}

The main result of \cite{AST79} is a beautiful class of explicit formulas for $P(\dg)$, non-commutative factorizations along normal orders of $\Delta(\dn^+)$.  For $\dsl_2$, it is an enjoyable exercise to prove
\begin{equation} \label{P sum}
   P(\dsl_2) = \sum_{k=0}^\infty {\ts\frac{1}{k!}} (-1)^k F^k E^k\, \prod_{i=1}^k (H + 1 + i)^{-1}.
\end{equation}

The AST factors of $P(\dg)$ generalize $P(\dsl_2)$.  For $t \in \bC$, define
\begin{equation} \label{Pt sum}
   Q_t(\dsl_2) := \sum_{k=0}^\infty {\ts\frac{1}{k!}} (-1)^k 
   F^k E^k\, \prod_{i=1}^k (H + t + i)^{-1} \in \dF(\dsl_2).
\end{equation}
Note that $Q_1(\dsl_2) = P(\dsl_2)$.  For $\alpha \in \Delta(\dn^+)$, let $\da_\alpha$ be the copy of $\dsl_2$ in $\dg$ corresponding to $\alpha$, and let $\dl_\alpha$ be its augmentation by $\dh$:
\begin{equation} \label{a and l}
   \da_\alpha := \Span_\bC \bigl\{ E_\alpha, F_\alpha, H_\alpha \bigr\}, \qquad
   \dl_\alpha := \dh + \da_\alpha.
\end{equation}
We remark that $P(\dl_\alpha)$ and $P(\da_\alpha)$ are the same.  Indeed, $P(\dg) = P(\dg_\ss)$ by definition.

\meno {\bf Definition.}
Fix a normal order $(\alpha_1, \ldots, \alpha_m)$ of $\Delta(\dn^+)$.  For $\tau \in \dh^*$, set
\begin{equation*}
   Q_\tau(\dg) := \prod_{r=1}^m Q_{\tau (H_{\alpha_r})} (\da_{\alpha_r}).
\end{equation*}

\begin{thm} \cite{AST79} \label{fin fac}
For any normal ordering of\/ $\Delta(\dn^+)$, $P(\dg) = Q_{\rho_\dg}(\dg)$.
\end{thm}

Observe that $\rho_\dg(H_\alpha) \in \bZ^+$ for $\alpha \in \Delta(\dn^+)$, so those $Q_t$ occurring as AST factors of $P(\dg)$ have $t \in \bZ^+$.  Theorem~\ref{fin fac} is in fact a corollary of the following more general result, discovered later by Zhelobenko.

\begin{thm} \cite{Zh90}
For all $\tau \in \dh^*$, $Q_\tau(\dg)$ is independent of the choice of normal order of\/ $\Delta(\dn^+)$.
\end{thm}

\subsection{Denominators} \label{Denoms}

In an obvious sense, the {\em total denominator\/} of each of the factorizations of $P(\dg)$ in Theorem~\ref{fin fac} is the commutative formal product
\begin{equation} \label{abs denom}
   D(\dg) := \prod_{i=1}^\infty \, \prod_{\alpha \in \Delta(\dn^+)} \, 
   \bigl(H_\alpha + \rho_\dg (H_\alpha) + i \bigr).
\end{equation}
This has the following implication: if $V$ is any \r\ in $\O(\dg)$ and $\mu$ is any weight on which no factor of $D(\dg)$ is zero, then the formula of Theorem~\ref{fin fac} defines an action of $P(\dg)$ on $V_\mu$.

\begin{prop} \label{abs denom prop}
$D(\dg)$ divides the total denominator of any formula for $P(\dg)$.
\end{prop}

\meno {\em Idea of proof.\/}
Suppose that $\mu$ is a weight annihilating some factor $\bigl( H_\alpha + \rho_\dg (H_\alpha) + i \bigr)$ of $D(\dg)$.  By~(\ref{Pg domain}), it suffices to find an object $V$ of $\O(\dg)$ with $V_\mu^{\dn^+} \cap \dn^- V \not= 0$.

Let $\cdot$ denote the affine {\em dot action\/}
\begin{equation*}
   w \cdot \mu := w (\mu + \rho_\dg) - \rho_\dg
\end{equation*}
of $W(\dg)$ on $\dh^*$.  Then $s_\alpha \cdot \mu - \mu = i \alpha$, so by a well known result of Bernstein, Gel'fand, and Gel'fand, the Verma module $M(\dg, s_\alpha \cdot \mu)$ of $\dg$ with highest weight $s_\alpha \cdot \mu$ satisfies
\begin{equation*}
   M(\dg, s_\alpha \cdot \mu)_\mu^{\dn^+} \cap \dn^- M(\dg, s_\alpha \cdot \mu) \not= 0. \qquad \Box
\end{equation*}

\subsection{Infinite commutative factorizations} \label{ICFs}

Let $\dZ(\dg)$ be the center of $\dU(\dg)$.  Zhelobenko discovered an infinite commutative factorization of $P(\dg)$ built from the Casimir element $\Omega_2$ of $\dZ(\dg)$.  In order to describe it we must extend the dot action of $W(\dg)$ to $\dU(\dh)$ and define the {\em shift action\/} of $\dh^*$ on $\dU(\dh)$.  Regard elements $\Uhelt$ of $\dU(\dh)$ as polynomials on $\dh^*$, and for $w \in W(\dg)$ and $\nu, \mu \in \dh^*$, set
\begin{equation} \label{shift}
   (w \Uhelt)(\mu) := \Uhelt (w^{-1} \mu), \quad
   (w \cdot \Uhelt)(\mu) := \Uhelt (w^{-1} \cdot \mu), \quad
   \Uhelt^\nu(\mu) := \Uhelt (\nu + \mu).
\end{equation}
Write $\dU(\dh)^{W(\dg) \cdot}$ for the subalgebra of $\dU(\dh)$ invariant under the dot action.

Recall that the subalgebra $\dU(\dg)^\dh = \dU(\dg)_0$ of $\dU(\dg)$ decomposes as $\bigl( \dn^- \dU(\dg) \dn^+ \bigr)_0 \oplus \dU(\dh)$, where $\bigl( \dn^- \dU(\dg) \dn^+ \bigr)_0$ is a two-sided ideal.  The {\em Harish-Chandra homomorphism\/} $\HC_\dg$ is the associated projection from $\dU(\dg)^\dh$ to $\dU(\dh)$.  By a well known result of Harish-Chandra, it restricts to an isomorphism
\begin{equation*}
   \HC_\dg: \dZ(\dg) \to \dU(\dh)^{W(\dg) \cdot}.
\end{equation*}

\begin{thm} \cite{Zh93} \label{inf comm fac}
As an element of\/ $\End_{\rho(\dh)} M(\dg)$,
\begin{equation*}
   P(\dg)\ = \prod_{\nu \in \Delta(\dU^+(\dn^+))} 
   \frac {\Omega_2 - (\HC_\dg \Omega_2)^\nu} {(\HC_\dg \Omega_2) - (\HC_\dg \Omega_2)^\nu}\,.
\end{equation*}
\end{thm}

\meno {\em Idea of proof.\/}
The infinite product is interpreted as follows: if the factors are applied successively to any $\Theta \in M(\dg)$, the resulting sequence eventually stabilizes.  To prove that it stabilizes at $P(\dg) \Theta$, recall~(\ref{h decomp of M}) and check that $\Omega_2$ acts on $M(\dg)_{-\nu}$ by $(\HC_\dg \Omega_2)^\nu$ for all $\nu \in \dU(\dn^+)$.  Therefore the $\nu$-factor of the infinite product acts by $0$ on $M(\dg)_{-\nu}$ and by $1$ on $M(\dg)_0$, so the entire product acts by $1$ on $M(\dg)_0$ and by $0$ on all other weight spaces.  Now apply~(\ref{Herm decomp of M}) and Theorem~\ref{Verma}.  $\Box$

\medbreak
It was observed in \cite{CS05} that for $\dg$ simple, Theorem~\ref{inf comm fac} holds for any non-constant element $\Omega$ of $\dZ(\dg)$ replacing $\Omega_2$.  For $\dsl_2$, it may be rewritten as
\begin{equation*}
   P(\dsl_2) = \prod_{i=1}^\infty \Bigl( 1 - \frac{FE} {i (H + 1 + i)} \Bigr).
\end{equation*}

It is an intriguing fact that the AST factors $Q_t$ occurring in Theorem~\ref{fin fac} have themselves a similar infinite commutative factorization.  In a natural telescopic sense explained in Theorem~15 of \cite{CS05}, for $t \in \bZ^+$ we have
\begin{equation} \label{Pt fac}
   Q_t(\dsl_2) = \prod_{i=t}^\infty \Bigl( 1 - \frac{F E} {i (H + 1 + i)} \Bigr).
\end{equation}
In particular, $Q_t(\dsl_2)$ annihilates all but the highest~$t$ weight spaces of $M(\dsl_2)$: its image is $\bigoplus_{i = 0}^{t - 1} M(\dsl_2)_{-2i}$.  However, for $t > 1$ it is not the Hermitian projection operator onto this sum.

%%%%%%%%%%%%%%%%
\section{The relative extremal projector}  \label{REP} 
%%%%%%%%%%%%%%%%

Let $\dl$ be a {\em standard\/} reductive subalgebra of $\dg$, \ie\ the Levi subalgebra of a standard parabolic subalgebra.  Thus $\dl$ contains $\dh$ and has triangular decomposition $\dl^- \oplus \dh \oplus \dl^+$, where $\dl^\pm := \dl \cap \dn^\pm$, and its positive root system $\Delta(\dl^+)$ has simple system
\begin{equation*}
   \Pi(\dl^+) = \Delta(\dl^+) \cap \Pi(\dn^+).
\end{equation*}
Let $\du^- \oplus \dl \oplus \du^+$ be the $\dl$-invariant decomposition of $\dg$ such that $\du^\pm \subset \dn^\pm$.

The relative \expro\ is introduced in the next two theorems, which are parallel to Theorems~\ref{expro dfn} and~\ref{Verma}.

\begin{thm} \label{rel expro dfn} \cite{CS03}
There is a unique non-zero Hermitian idempotent $P(\dg, \dl)$ in $\dF(\dg)_0$, the\/ {\em relative \expro,} which commutes with $\dl$ and satisfies
\begin{equation*}
   \du^+ P(\dg, \dl) = 0 = P(\dg, \dl) \du^-.
\end{equation*}
\end{thm}

Recall that $\dF(\dg)$ is identified with $\End_{\rho(\dh)} M(\dg)$.  The relative analog of~(\ref{Herm decomp of M}) was observed in Lemma~1 of \cite{CS03}.  In order to state it, we must explain how $M(\dg)$ decomposes as a direct sum of copies of $M(\dl)$.

Note that those monomials $F^I$ from~(\ref{roots}) which are contained in $\dU(\du^-)$ form a basis of it.  The \expro\ $P(\dl)$ of $\dl$ maps each of them to a non-zero $\dl$-\hwv\ $P(\dl) (F^I)$ in $M(\dg)$.  Applying $\dU(\dl^-) \dF(\dh)$ to this \hwv\ gives an $\dl$-submodule of $M(\dg)$ isomorphic to $M(\dl)$, and $M(\dg)$ is the Hermitian direct sum of these submodules:
\begin{equation} \label{l decomp of M}
   M(\dg) = \bigoplus_{F^I \in \dU(\du^-)} \dU(\dl^-) \dF(\dh) P(\dl) (F^I).
\end{equation}

The {\em highest $\dl$-submodule\/} $M(\dg)^{\du^+}$ of $M(\dg)$ is the copy of $M(\dl)$ generated by~$1$, and the lower submodules are the other copies:
\begin{align} \label{rel Herm decomp of M}
   & M(\dg)^{\du^+} = \dU(\dl^-) \dF(\dh), \nonumber \\[6pt]
   & \du^- M(\dg) = \bigoplus_{F^I \in \dU^+(\du^-)} \dU(\dl^-) \dF(\dh) P(\dl) (F^I), \\[6pt]
   & M(\dg) = M(\dg)^{\du^+} \oplus \du^- M(\dg). \nonumber
\end{align}

\begin{thm} \label{rel Verma} \cite{CS03}
As an element of\/ $\End_{\rho(\dh)} M(\dg)$, the relative \expro\ $P(\dg, \dl)$ is the projection from $M(\dg)$ to $M(\dg)^{\du^+}$ with kernel\/ $\du^- M(\dg)$.
\end{thm}

Because $P(\dg, \dl)$ commutes with $\dl$, we seek formulas for it whose terms lie in $\dF(\dg)^\dl$, the commutant of $\dl$ in $\dF(\dg)$.  In the relative case, the natural analogs of the factorizations of $P(\dg)$ discussed above have factors with numerators in $\dU(\dg)^\dl$ and denominators in the center $\dZ(\dl)$ of $\dU(\dl)$.  Such denominators are allowed because by Lemma~2 of \cite{CS03}, non-zero elements of $\dZ(\dl)$ are invertible in $\dF(\dg)$.

\meno {\bf Remark.}
By Theorem~6 of \cite{CS03}, if $\dl'$ is a standard reductive subalgebra of $\dl$ then
\begin{equation} \label{successive}
   P(\dg, \dl') = P(\dg, \dl) P(\dl, \dl') = P(\dl, \dl') P(\dg, \dl).
\end{equation}
In particular, $\dl' = \dh$ gives $P(\dg) = P(\dg, \dl) P(\dl) = P(\dl) P(\dg, \dl)$, as $P(\dg, \dh) = P(\dg)$.

Compare this factorization of $P(\dg)$ to the AST factorizations from Theorem~\ref{fin fac}.  Any normal order of $\Delta(\dl^+)$ can be extended to a normal order of $\Delta(\dn^+)$ with $\Delta(\dl^+)$ all to the left or all to the right.  For such normal orders of $\Delta(\dn^+)$, the product of those AST factors $Q_{\rho_\dg(H_\alpha)}(\da_\alpha)$ of $P(\dg)$ with $\alpha \in \Delta(\dl^+)$ is equal to $P(\dl)$.  However, one cannot ``cancel'' this factor $P(\dl)$ of $P(\dg)$ with the one in~(\ref{successive}) for $\dl' = \dh$: $P(\dg, \dl)$ is not in general the product of the non-$\Delta(\dl^+)$ AST factors of $P(\dg)$.

For example, in obvious notation, let $\{\alpha_{12}, \alpha_{13}, \alpha_{23} \}$ be a normal order of the positive roots of $\dsl_3$, and write $\da_{ij}$ and $\dl_{ij}$ for the subalgebras in~(\ref{a and l}).  Then~(\ref{successive}) and Theorem~\ref{fin fac} give
\begin{equation*}
   P(\dsl_3) = P(\dsl_3, \dl_{23}) P(\dl_{23}) = P(\dl_{12}) Q_2(\da_{13}) P(\dl_{23}),
\end{equation*}
but $P(\dsl_3, \dl_{23})$ is not equal to $P(\dl_{12}) Q_2(\da_{13})$.  Indeed, $P(\dsl_3, \dl_{23})$ projects $M(\dsl_3)$ to $\bigoplus_{n=0}^\infty M(\dsl_3)_{-n \alpha_{23}}$, so it annihilates $F_{\alpha_{13}}$.  On the other hand, $P(\dl_{12}) Q_2(\da_{13})$ does not: $Q_2(\da_{13}) (F_{\alpha_{13}})$ is a non-zero $\dF(\dh)$-multiple of $F_{\alpha_{13}}$, and $P(\dl_{12}) (F_{\alpha_{13}}) \not= 0$.

The following lemma is in a sense a of converse of~(\ref{successive}).

\begin{lemma} \label{PglPl}
$P(\dg, \dl)$ is the unique element of\/ $\dF(\dg)^\dl$ such that $P(\dg, \dl) P(\dl) = P(\dg)$.
\end{lemma}

\meno {\em Proof.\/}
We saw in~(\ref{l decomp of M}) that as an $\dl$-module, $M(\dg)$ is a direct sum of copies of $M(\dl)$.  Each copy of $M(\dl)$ is generated under $\dl$ by its $\dl$-\hwv, and $P(\dl) M(\dg)$ is the space of all such \hwv s.  Thus any element $\pi$ of $\dF(\dg)^\dl$ is determined by its action on $P(\dl) M(\dg)$.  In particular, if $\pi P(\dl) = P(\dg)$, then $\pi$ must be $P(\dg, \dl)$ because $P(\dg, \dl) P(\dl) = P(\dg)$.  $\Box$

\subsection{Infinite commutative factorizations}

In \cite{CS03} and \cite{CS05} we give infinite commutative factorizations of $P(\dg, \dl)$.  Theorem~7 of \cite{CS03} is the relative analog of Theorem~\ref{inf comm fac}, a factorization built from the Casimir element $\Omega_2$.  Theorem~4~(3) of \cite{CS05} shows that it holds with almost any element of $\dZ(\dg)$ replacing $\Omega_2$.

Following Section~\ref{ICFs}, write $\dU(\dh)^{W(\dl) \cdot}$ for the subalgebra of $\dU(\dh)$ invariant under the $\dl$-dot action of the Weyl group $W(\dl)$ of $\dl$.  The Harish-Chandra isomorphism $\HC_\dl: \dZ(\dl) \to \dU(\dh)^{W(\dl) \cdot}$ extends to an isomorphism from $\Frac \bigl( \dZ(\dl) \bigr)$ to $\dF(\dh)^{W(\dl) \cdot}$.  We extend it further by $t \mapsto t$ to an isomorphism
\begin{equation*}
   \HC_\dl: \Frac \bigl( \dZ(\dl) \bigr) [t] \to \dF(\dh)^{W(\dl) \cdot} [t].
\end{equation*}

It is an elementary but crucial observation that although in general $\rho_\dl \not= \rho_\dg$, the $\dl$- and $\dg$-dot actions of $W(\dl)$ are the same, because $W(\dl)$ stabilizes $\rho_\dg - \rho_\dl$.  Therefore we may speak unambiguously of the dot action of $\dl$.  Taking $w \in W(\dg)$, $\nu \in \dh^*$, and $\Uhelt \in \dF(\dh)$, let us observe that
\begin{equation} \label{wQnu}
   w \cdot (\Uhelt^\nu) = (w \cdot \Uhelt)^{w \nu} = (w\Uhelt)^{w \nu - w \rho_\dg + \rho_\dg}, \qquad
   \dU(\dh)^{W(\dg) \cdot} \subseteq \dU(\dh)^{W(\dl) \cdot}.
\end{equation}

\begin{thm} \label{any Omega} \cite{CS05}
Let $\Omega$ be any element of\/ $\dZ(\dg)$ that has a non-constant component over every simple summand of\/ $\dg$.  Then as elements of\/ $\End_{\rho(\dh)} M(\dg)$,
\begin{equation*}
   P(\dg, \dl) = \biggl( \HC_\dl^{-1} \prod_{\nu \in \Delta(\dU^+(\du^+))} 
   \frac {t - (\HC_\dg \Omega)^\nu} {(\HC_\dg \Omega) - (\HC_\dg \Omega)^\nu}
   \biggr) \bigg|_{t = \Omega} \,.
\end{equation*}
\end{thm}

\meno {\em Idea of proof.\/}
We first explain the expression.  $\HC_\dg \Omega$ is dot-invariant, so (\ref{wQnu}) gives
\begin{equation*}
   w \cdot \biggl( \frac {t - (\HC_\dg \Omega)^\nu} {(\HC_\dg \Omega) - (\HC_\dg \Omega)^\nu} \biggr)
   =    \frac {t - (\HC_\dg \Omega)^{w \nu}} {(\HC_\dg \Omega) - (\HC_\dg \Omega)^{w \nu}} \,.
\end{equation*}
Since $W(\dl)$ leaves $\Delta \bigl(\dU^+(\du^+) \bigr)$ invariant and partitions it into finite orbits, the product may be written as an infinite product of finite products over these orbits.  Each finite product is $W(\dl)$ dot-invariant, so $\HC_\dl^{-1}$ may be applied to it to give an element of $\Frac \bigl( \dZ(\dl) \bigr) [t]$.  Then substituting $\Omega$ for $t$ gives an element of $\dZ(\dg) \Frac \dZ(\dl)$.  The product of all these elements is the right hand side.

The proof that this infinite product is $P(\dg, \dl)$ is similar to the proof of Theorem~\ref{inf comm fac}.  Its factors commute with $\dl$, so by~(\ref{l decomp of M}) and~(\ref{rel Herm decomp of M}) it suffices to prove that it acts by $1$ on $1$ and by $0$ on $P(\dl) (F^I)$ for all $F^I$ in $\dU^+(\du^-)$.  It is not hard to check that on $P(\dl) (F^I)$, $\Omega$ acts by $(\HC_\dg \Omega)^{|I|}$ and the elements of $\dZ(\dl)$ act by their images under $\HC_\dl$, so the $\nu$-factor contributes action
\begin{equation*}
   \frac {(\HC_\dg \Omega)^{|I|} - (\HC_\dg \Omega)^\nu} {(\HC_\dg \Omega) - (\HC_\dg \Omega)^\nu} \,.
\end{equation*}
This numerator is $0$ for $\nu = |I|$, so it only remains to prove that none of the denominators are~$0$.  This follows from the non-constancy condition on $\Omega$.  $\Box$

\medbreak
The total denominator $\HC_\dl^{-1} \prod_{\nu \in \Delta(\dU^+(\du^+))} \bigl( \HC_\dg \Omega - (\HC_\dg \Omega)^\nu \bigr)$ of the formula for $P(\dg, \dl)$ given in Theorem~\ref{any Omega} is a formal product of elements of $\dZ(\dl)$.  The formula is not efficient: its denominator is larger than necessary.  Theorem~3 of \cite{CS05} gives a general scheme for constructing infinite factorizations of $P(\dg, \dl)$ with factors in $\dZ(\dg) \Frac \dZ(\dl)$, as opposed to simply $\bC[\Omega] \Frac \dZ(\dl)$ for some $\Omega$ in $\dZ(\dg)$.  It is applied in Theorems~4~(1) and 4~(2) to give factorizations with smaller denominators.  These two parts of Theorem~4 are identical for $\dl$ maximal, the most important case in the context of finding a relative version of Theorem~\ref{fin fac}.  We now recall Theorem~4~(1).

The center $\dz(\dl)$ of $\dl$ is of course a subalgebra of $\dh$.  Define
\begin{equation*}
   \dz^+(\dl) := \bigl\{ T \in \dz(\dl): \mbox{\rm Real Part} \bigl(\alpha(T) \bigr) > 0
   \ \forall\ \alpha \in \Delta(\du^+) \bigr\}.
\end{equation*}

For $T \in \dh$, write $W(\dg)^T$ for the $W(\dg)$-stabilizer of $T$.  The stabilizer and dot-stabilizer of $T$ are the same, so the following polynomial is in $\dU(\dh)^{W(\dg) \cdot}[t]$:
\begin{equation*}
   p_T(t) := \prod_{w \in W(\dg) / W(\dg)^T} (t - w \cdot T).
\end{equation*}

\begin{thm} \label{4(1)} \cite{CS05}
For all $T \in \dz^+(\dl)$,
\begin{equation} \label{4(1) fac}
   P(\dg, \dl) = \prod_{c \in \Delta(\dU^+(\du^+)) (T)}
   \frac { \bigl( \HC_\dg^{-1} p_T(t) \bigr) \big|_{t = T + c} }
   { \HC_\dl^{-1} p_T(T + c) } \,.
\end{equation}
\end{thm}

\meno {\em Idea of proof.\/}
As in the proof of Theorem~\ref{any Omega}, the factors commute with $\dl$, so it suffices to prove that the product acts by $1$ on $1$ and by $0$ on $P(\dl) (F^I)$ for all $F^I$ in $\dU^+(\du^-)$.  Check that the action of the numerators and denominators on $P(\dl) (F^I)$ is multiplication by the following quantities:
\begin{align*}
   & \bigl( \HC_\dg^{-1} p_T(t) \bigr) \big|_{t = T + c} \rightsquigarrow
   p_T^{|I|}(T + c) = \prod_{W(\dg) / W(\dg)^T} \bigl(T + c - (w \cdot T)^{|I|} \bigr), \\[6pt]
   & \HC_\dl^{-1} p_T(T + c) \rightsquigarrow
   p_T(T + c) = \prod_{w \in W(\dg) / W(\dg)^T} (T + c - w \cdot T).
\end{align*}
These multipliers are equal for $I = 0$, and $T \in \dz^+(\dl)$ implies that the denominators never act by zero.  For $I >0$, the factor of the numerator's multiplier with $c = |I| (T)$ and $w = e$ acts by zero.  $\Box$

\subsection{Denominators}

As discussed in Section~\ref{Denoms} for $P(\dg)$, formulas for $P(\dg, \dl)$ with smaller denominators are better, as they have larger domains of definition.  We now recall Theorem~8 and Conjecture~1 of \cite{CS05}.  The theorem gives a lower bound for the denominator of $P(\dg, \dl)$.  It generalizes Proposition~\ref{abs denom prop}, and its proof is again a BGG argument.  The conjecture predicts that the lower bound can be achieved.

\begin{prop} \label{rel BGG}
The total denominator of any formula for $P(\dg, \dl)$ is divisible by
\begin{equation*}
   D(\dg, \dl) := \prod_{i=1}^\infty \, \HC_\dl^{-1} \, \biggl( \, \prod_{\alpha \in \Delta(\du^+)} \,
   \bigl(H_\alpha + i \bigr)^{\rho_\dg} \biggr).
\end{equation*}
\end{prop}

\begin{conj} \label{rel conj}
There is a formula for $P(\dg, \dl)$ with total denominator $D(\dg, \dl)$.
\end{conj}

For $\dl = \dh$, this follows from Theorem~\ref{fin fac}.  In Theorem~13 of \cite{CS05} we use Theorem~\ref{4(1)} to prove it in a few additional cases:

\begin{thm} \label{13}
Conjecture~\ref{rel conj} holds if\/ $\dg$ is of type $A_n$ or $B_n$ and the simple roots of\/ $\dl$ form a ``ray'': a connected segment of the Dynkin diagram of\/ $\dg$ including an end root, short in the case of $B_n$.
\end{thm}

\meno {\em Idea of proof.\/}
Check that $T - w \cdot T = (T - wT)^{\rho_\dg}$, so the denominator of~(\ref{4(1) fac}) is
\begin{equation*}
   D(\dg, \dl, T) := \prod_{c \in \Delta(\dU^+(\du^+)) (T)} \, \HC_\dl^{-1} \,
   \biggl( \prod_{W(\dg) / W(\dg)^T} \bigl( T - wT + c \bigr)^{\rho_\dg} \biggr).
\end{equation*}
If $\dl$ is maximal, $\du^+$ is \irr\ under $\dl$, and $|\Delta(\du^+)| = |W(\dg) / W(\dl)| - 1$, then $D(\dg, \dl, T)$ is proportional to $D(\dg, \dl)$ for any $T \in \dz^+(\dl)$.  However, these conditions hold \iff\ $\dg$ is of type $A_n$ or $B_n$ and the simple root missing from $\dl$ is an end root, long in the case of $B_n$.  An inductive argument based on Theorem~5 of \cite{CS05} now gives the result.  $\Box$

\meno {\bf Remark.}
In order to understand some of the obstacles to further progress, it may be helpful to consider the case that $\dg$ is $\dgo_5$ and $\dl$ is a long $\dgl_2$.  Let $\pm \ep_1$, $\pm \ep_2$, and $\pm \ep_1 \pm \ep_2$ be the roots of $\dgo_5$.  Take simple roots $\ep_1 - \ep_2$ and $\ep_2$ and let $\dl$ have simple root $\ep_1 - \ep_2$.  We may identify $\dh$ and $\dh^*$ via 
\begin{equation*}
   H_{\pm \ep_i} \equiv \pm 2\ep_i, \qquad
   H_{\pm \ep_i \pm \ep_j} \equiv \pm \ep_i \pm \ep_j.
\end{equation*}

When $\dl$ is maximal, $\dz(\dl)$ is 1-dimensional, so there is essentially only one choice of $T$ in Theorem~\ref{4(1)}.  Here that choice is $H_{\ep_1 + \ep_2}$.  Up to proportionality, the ``extra factors'' of $D(\dgo_5, \dl, H_{\ep_1 + \ep_2})$, those not occurring in $D(\dgo_5, \dl)$, are
\begin{equation*}
   \bigl( H_{\ep_1 + \ep_2} + \oh \bigr)^{\rho_\dg},\
   \bigl( H_{\ep_1 + \ep_2} + {\ts \frac{3}{2}} \bigr)^{\rho_\dg},\
   \bigl( H_{\ep_1 + \ep_2} + {\ts \frac{5}{2}} \bigr)^{\rho_\dg}, \ldots .
\end{equation*}
We expect that in this example, no factorization of $P(\dgo_5, \dl)$ with factors drawn from $\dZ(\dgo_5) \Frac \dZ(\dl)$ achieves the minimal denominator $D(\dgo_5, \dl)$; our guess is that it can only be attained by factorizations over $\dU(\dgo_5)^\dl \Frac \dZ(\dl)$.

\subsection{Non-commutative factorizations}

We now give some new results in the case that $\dg$ is $\dsl_4$ or $\dsl_5$: for any standard reductive subalgebra $\dl$ of $\dsl_4$ and for all but one such subalgebra of $\dsl_5$, we prove that $P(\dg, \dl)$ has non-commutative factorizations analogous to the factorizations of $P(\dg)$ given in Theorem~\ref{fin fac}.  In these factorizations, the factors of $P(\dg, \dl)$ are indexed by certain reductive subalgebras $\dm$ of $\dg$, not in general standard in $\dg$, which themselves contain $\dl$ as a maximal standard reductive subalgebra.  The factor $Q(\dm, \dl)$ corresponding to $\dm$ is an element of $\dF(\dm)^\dl$.

Unlike the method of \cite{AST79}, the method we will present is non-constructive: it only shows that the factors $Q(\dm, \dl)$ exist.  They are relative analogs of the AST factors $Q_t(\da_\alpha)$, and it would be interesting to have explicit formulas for them such as~(\ref{Pt sum}) and~(\ref{Pt fac}).  Formulas for their total denominators as formal products in $\dZ(\dl)$ would allow a resolution of Conjecture~\ref{rel conj}.

If $\dl$ is a maximal standard subalgebra of $\dg$, our method gives no non-trivial factorizations of $P(\dg, \dl)$.  At the other extreme, for $\dl = \dh$ it gives only a weaker version of Theorem~\ref{fin fac}.  Thus the interest lies in the cases $1 \le |\Pi(\dl^+)| \le |\Pi(\dn^+)| - 2$, that is, $1 \le \rank(\dl_\ss) \le \rank(\dg_\ss) - 2$.

Let us recall some standard notation for $\dsl_n$.  Take $\dn^+$ and $\dn^-$ to be the upper and lower triangular matrices, respectively, and $\dh$ to be the diagonal matrices.  Writing $e_{ij}$ for the usual elementary $n \times n$ matrix, the positive and negative root vectors and corresponding elements of $\dh$ are
\begin{equation*}
   E_{ij} := e_{ij}, \quad F_{ij} := e_{ji}, \quad H_{ij} := e_{ii} - e_{jj}; \quad 1 \le i < j \le n.
\end{equation*}
Let $\ep_i$ be the $i^\thup$ standard basis vector of $\bR^n$, so that $E_{ij}$ has root $\alpha_{ij} := \ep_i - \ep_j$.

For $1 \le i_1 < i_2 < \cdots < i_r \le n$, define subalgebras $\da_{i_1 \cdots i_r}$ and $\dl_{i_1 \cdots i_r}$ of $\dsl_n$ by
\begin{equation*}
   \da_{i_1 \cdots i_r} := \Span_\bC
   \bigl\{ E_{i_a i_b},\, F_{i_a i_b},\, H_{i_a i_b}: 1 \le a < b \le r \bigr\}, \quad
   \dl_{i_1 \cdots i_r} := \dh + \da_{i_1 \cdots i_r}.
\end{equation*}
Note that $\dl_{i_1 \cdots i_r}$ is standard \iff\ $i_1, \ldots, i_r$ are consecutive.

If $\{ j_1 < \cdots < j_s \} \subseteq \{i_1 < \cdots < i_r \}$, then $\dl_{j_1 \cdots j_s} \subseteq \dl_{i_1 \cdots i_r}$.  In this case we write
\begin{equation*}
   P_{i_1 \cdots i_r}^{j_1 \cdots j_s} := P(\dl_{i_1 \cdots i_r}, \dl_{j_1 \cdots j_s}), \qquad
   P_{1 \cdots n}^{ j_1 \cdots j_s} := P(\dsl_n, \dl_{j_1 \cdots j_s}).
\end{equation*}

At the other extreme, if $\{i_1 < \cdots < i_r \}$ and $\{ j_1 < \cdots < j_s \}$ are disjoint, then $\dl_{i_1 \cdots i_r}$ and $\dl_{j_1 \cdots j_s}$ commute.  In this case $\dsl_n$ has the reductive subalgebra 
\begin{equation*}
   \dl_{i_1 \cdots i_r, j_1 \cdots j_s} := \dl_{i_1 \cdots i_r} + \dl_{j_1 \cdots j_s},
\end{equation*}
and $P(\dl_{i_1 \cdots i_r, j_1 \cdots j_s}, \dl_{j_1 \cdots j_s})$ is simply $P_{i_1 \cdots i_r}$.

We now state our results; their proofs are given in Section~\ref{Proofs}.  Keep in mind that in these factorizations of $P(\dg, \dl)$, the factors commute with $\dl$ but not always with each other.  Some of them do commute, and the reader will note that their possible orders are closely related to normal orders of $\Delta(\dn^+)$.  Some of the factors coincide with the AST operators $Q_t(\dsl_2)$ in~(\ref{Pt sum}), and so for $1 \le a < b \le n$ we define
\begin{equation*}
   Q_{a b} := Q_{b - a}(\da_{a b}).
\end{equation*}

For $\dsl_4$ we are concerned only with the case $|\Pi(\dl^+)| = 1$, and for $\dsl_5$ only with the cases $|\Pi(\dl^+)| = 1$ or~$2$.  For both $\dsl_4$ and $\dsl_5$, up to isomorphism the only choices of $\dl$ with $|\Pi(\dl^+)| = 1$ are $\dl_{12}$ and $\dl_{23}$.  For $\dsl_5$ there are four choices with $|\Pi(\dl^+)| = 2$:
\begin{equation*}
   \dl_{123}, \qquad \dl_{234}, \qquad \dl_{12, 34}, \qquad
   \dl_{12, 45}.
\end{equation*}
We will use obvious notation such as $P_{12345}^{12, 34}$ for $P(\dsl_5, \dl_{12, 34})$.  The one case we will not treat is that of $P_{12345}^{12, 45}$; the reason for this is explained in Section~\ref{Remarks}.

\begin{thm} \label{sl4}
For $|\Pi(\dl^+)| = 1$, $P(\dsl_4, \dl)$ has the following factorizations:

\begin{enumerate}

\item[(i)]
For $\dl = \dl_{12}$, there is a unique element $Q_{124}^{12}$ of\/ $\dF(\dl_{124})^{\dl_{12}}$ such that
\begin{equation*}
   P_{1234}^{12} = P_{123}^{12}\, Q_{124}^{12}\, P_{34}.
\end{equation*}

\smallbreak \item[(ii)]
For $\dl = \dl_{23}$, $Q_{14}$ is the unique element of\/ $\dF(\dl_{14})^\dh$ such that
\begin{equation*}
   P_{1234}^{23} = P_{123}^{23}\, Q_{14}\, P_{234}^{23}.
\end{equation*}

\end{enumerate}
\end{thm}

\begin{thm} \label{sl5r1}
For $|\Pi(\dl^+)| = 1$, $P(\dsl_5, \dl)$ has the following factorizations:

\begin{enumerate}

\item[(i)]
For $\dl = \dl_{12}$, there is a unique element $Q_{125}^{12}$ of\/ $\dF(\dl_{125})^{\dl_{12}}$ such that
\begin{equation*}
   P_{12345}^{12} = P_{1234}^{12}\, Q_{125}^{12}\, P_{345} =
   P_{123}^{12}\, Q_{124}^{12}\, P_{34}\, Q_{125}^{12}\, Q_{35}\, P_{45}.
\end{equation*}

\smallbreak \item[(ii)]
Let $Q_{235}^{23} \in \dF(\dl_{235})^{\dl_{23}}$ be $Q_{124}^{12}$ with all indices shifted up one.  For $\dl = \dl_{23}$, $Q_{15}$ is the unique element of\/ $\dF(\dl_{15})^{\dh}$ such that
\begin{equation*}
   P_{12345}^{23} = P_{1234}^{23}\, Q_{15}\, P_{2345}^{23} =
   P_{123}^{23}\, Q_{14}\, P_{234}^{23}\, Q_{15}\, Q_{235}^{23}\, P_{45}.
\end{equation*}

\end{enumerate}
\end{thm}

\begin{thm} \label{sl5r2}
For $|\Pi(\dl^+)| = 2$, $P(\dsl_5, \dl)$ has the following factorizations:

\begin{enumerate}

\item[(i)]
For $\dl = \dl_{123}$, there is a unique element $Q_{1235}^{123}$ of\/ $\dF(\dl_{1235})^{\dl_{123}}$ such that
\begin{equation*}
   P_{12345}^{123} = P_{1234}^{123}\, Q_{1235}^{123}\, P_{45}.
\end{equation*}

\smallbreak \item[(ii)]
For $\dl = \dl_{234}$, $Q_{15}$ is the unique element of\/ $\dF(\dl_{15})^\dh$ such that
\begin{equation*}
   P_{12345}^{234} = P_{1234}^{234}\, Q_{15}\, P_{2345}^{234}.
\end{equation*}

\smallbreak \item[(iii)]
For $\dl = \dl_{12, 34}$, $Q_{125}^{12}$ above is the unique element of\/ $\dF(\dl_{125})^{\dl_{12}}$ such that
\begin{equation*}
   P_{12345}^{12, 34} = P_{1234}^{12, 34}\, Q_{125}^{12}\, P_{345}^{34}.
\end{equation*}

\end{enumerate}
\end{thm}

\section{Proofs} \label{Proofs}

In this section we prove Theorems~\ref{sl4}, \ref{sl5r1}, and~\ref{sl5r2}.  We will need the generalization of Lemma~1 of \cite{CS03}, given above as~(\ref{l decomp of M}), to reductive subalgebras $\dm$ of $\dg$ which contain $\dh$ but are not necessarily standard.  Such $\dm$ have triangular decomposition $\dm^- \oplus \dh \oplus \dm^+$, where $\dm^\pm := \dm \cap \dn^\pm$, but $\Delta(\dm^+) \cap \Pi(\dn^+)$ is not necessarily a simple system of the positive system $\Delta(\dm^+)$.

Throughout we will work over the field $\dF(\dh)$, and we abbreviate the phrase ``\hwv'' to HWV.  Remembering~(\ref{roots}), define
\begin{equation*}
   \I(\dg, \dm) := \bigr\{ I \in \bN^m:\, I_r = 0\ \forall\ \alpha_r \in \Delta(\dm^+) \bigr\}.
\end{equation*}

\begin{lemma} \label{nonstandard}
The set $\bigl\{ P(\dm) (F^I): I \in \I(\dg, \dm) \bigr\}$ is an $\dF(\dh)$-basis of $M(\dg)^\dm$.  For $I \in \I(\dg, \dm)$, the space $\dU(\dm^-) \dF(\dh) P(\dm) (F^I)$ is an $\dm$-submodule of $M(\dg)$ isomorphic to $M(\dm)$.  Moreover, $M(\dg)$ is the direct sum of these submodules:
\begin{equation} \label{m decomp of M}
   M(\dg) = \bigoplus_{I \in \I(\dg, \dm)} \dU(\dm^-) \dF(\dh) P(\dm) (F^I).
\end{equation}
\end{lemma}

\meno {\em Proof.\/}
If $P(\dm) (F^I)$ is non-zero then it is an $\dm$-HWV, so $\dU(\dm^-) \dF(\dh) P(\dm) (F^I)$ is $\dm$-isomorphic to $M(\dm)$ because $\dn^-$ acts freely on $M(\dg)$.

A PBW argument shows that the weight space dimensions on the right side of~(\ref{m decomp of M}) are no bigger than those on the left, with equality only if the sum is direct.  Conversely, if the right side contains $F^I$ for all $I \in \I(\dg, \dm)$, then it is $M(\dg)$.  To finish, induct on the usual partial order on $\Delta \bigl( \dU(\dn^+) \bigr)$: check that $P(\dm) (F^I) \equiv F^I$ modulo the sum of those $\dU(\dm^-) \dF(\dh) P(\dm) (F^J)$ with $|J| < |I|$.  $\Box$

\meno {\bf Definition.}
For $I \in \I(\dg, \dm)$, let $P(\dg, \dm, F^I)$ be the projection of $M(\dg)$ to the copy $\dU(\dm^-) \dF(\dh) P(\dm) (F^I)$ of $M(\dm)$ along the other summands of~(\ref{m decomp of M}).  By Theorem~\ref{Verma} and the $\dm$-invariance of~(\ref{m decomp of M}), $P(\dg, \dm, F^I)$ is an element of $\dF(\dg)^\dm$.

\meno {\bf Remark.}
We have just seen that it is possible to define the relative \expro\ $P(\dg, \dm)$ even for $\dm$ non-standard: it is $P(\dg, \dm, 1)$.

\subsection{$\dsl_3$: Warm up exercise} \label{sl3 proof}

Before proving the theorems, we illustrate the strategy by showing that there exists a unique element $\t Q_{13}$ of $\dF(\dl_{13})^\dh$ such that the projector $P(\dsl_3) = P_{123}$ factors as $P_{12} \t Q_{13} P_{23}$.  Of course Theorem~\ref{fin fac} tells us that $\t Q_{13}$ exists and is $Q_{13}$, but it is useful begin in the simplest setting.  We break the argument into several steps, which will be mirrored in the proofs of the theorems.

\meno {\em Step~1.\/}
By~(\ref{l decomp of M}), all $\dl_{12}$-HWVs in $M(\dsl_3)$ have weights in $-\Span_\bN \{\alpha_{13}, \alpha_{23} \}$, and all $\dl_{23}$-HWVs in $M(\dsl_3)$ have weights in $-\Span_\bN \{\alpha_{12}, \alpha_{13} \}$.  It follows that for {\em any\/} $\t Q_{13}$ in $\dF(\dl_{13})_0$, $P_{12} \t Q_{13} P_{23}$ annihilates all weight spaces $M(\dsl_3)_{-\nu}$ with
\begin{equation*}
   \nu \not\in \Span_\bN \bigl\{\alpha_{13}, \alpha_{23} \bigr\} \cap 
   \Span_\bN \bigl\{ \alpha_{12}, \alpha_{13} \bigr\} = \bN \alpha_{13}.
\end{equation*}
Therefore it suffices to choose $\t Q_{13}$ so that $P_{12} \t Q_{13} P_{23}$ maps~$1$ to $1$ and annihilates $M(\dsl_3)_{-n \alpha_{13}}$ for $n > 0$.

\meno {\em Step~2.\/}
By the PBW theorem, $\{ F_{12}^j F_{23}^j F_{13}^{n-j} \}_j$ and $\{ F_{23}^j F_{12}^j F_{13}^{n-j} \}_j$ are both bases of $M(\dsl_3)_{-n \alpha_{13}}$.  Therefore 
\begin{equation} \label{sl3 step2}
   P_{12} \bigl( M(\dsl_3)_{-n \alpha_{13}} \bigr) = \dF(\dh) P_{12} (F_{13}^n), \quad
   P_{23} \bigl( M(\dsl_3)_{-n \alpha_{13}} \bigr) = \dF(\dh) P_{23} (F_{13}^n).
\end{equation}
By~(\ref{l decomp of M}), both of these spaces are non-zero.  By the second of the two equations, we are done if we prove that there is a unique choice of $\t Q_{13}$ such that $P_{12} \t Q_{13}$ maps $1$ to $1$ and $P_{23}(F_{13}^n)$ to $0$ for $n > 0$.

\meno {\em Step~3.\/}
By a weight argument, Lemma~\ref{nonstandard} implies that
\begin{equation*}
   \bigoplus_{j = 0}^\infty M(\dsl_3)_{-j \alpha_{13}} =
   \bigoplus_{j = 0}^\infty \dU(\dl_{13}^-) \dF(\dh) P_{13} (F_{12}^j F_{23}^j).
\end{equation*}
In particular, for unique elements $h_0, \ldots, h_n$ of $\dF(\dh)$,
\begin{equation} \label{sl3 step3}
   P_{23} (F_{13}^n) = \sum_{i=0}^n h_i F_{13}^i P_{13} (F_{12}^{n-i} F_{23}^{n-i}).
\end{equation}

\meno {\em Step~4.\/}
The operator $P(\dl_{13}, \dh, F_{13}^k)$ in $\dF(\dl_{13})$ defined after Lemma~\ref{nonstandard} projects $M(\dl_{13})$ to its weight space $\dF(\dh) F_{13}^k = M(\dl_{13})_{-k \alpha_{13}}$ along its other weight spaces.  Abbreviate it to $P_{13}[k]$.  For any $\t Q_{13}$ in $\dF(\dl_{13})_0$, there exist unique elements $q_k$ of $\dF(\dh)$ such that
\begin{equation} \label{sl3 step4}
   \t Q_{13} = \sum_{k=0}^\infty q_k P_{13}[k].
\end{equation}

It is a crucial point that $P_{13}[k] \bigl(F_{13}^i P_{13} (F_{12}^j F_{23}^j) \bigr) = \delta_{i, k} F_{13}^i P_{13} (F_{12}^j F_{23}^j)$, because elements of $\dF(\dl_{13})$ such as $P_{13}[k]$ see $\dl_{13}$-HWVs such as $P_{13} (F_{12}^j F_{23}^j)$ as~$1$.  Therefore
\begin{equation*}
   \t Q_{13} P_{23} (F_{13}^n) = \sum_{i=0}^n q_i h_i F_{13}^i P_{13} (F_{12}^{n-i} F_{23}^{n-i}).
\end{equation*}

\meno {\em Step~5.\/}
Now apply $P_{12}$: we must choose $q_0, q_1, q_2, \ldots$ so that
\begin{equation} \label{sl3 step5}
   P_{12} \t Q_{13} P_{23} (F_{13}^n) = \sum_{i=0}^n
   q_i h_i P_{12} \bigl( F_{13}^i P_{13} (F_{12}^{n-i} F_{23}^{n-i}) \bigr)
\end{equation}
is~$1$ for $n=0$ and~$0$ for all $n > 0$.  We choose them successively.  Clearly $q_0$ must be~$1$.  Suppose that $q_1, \ldots, q_{n-1}$ have been determined.  By the first equation in~(\ref{sl3 step2}), every summand of~(\ref{sl3 step5}) is a multiple of $P_{12} (F_{13}^n)$.  Therefore there is a unique choice of $q_n$ such that the right side is zero \iff\ the coefficient of $q_n$ is non-zero.  Since $P_{12} (F_{13}^n) \not= 0$, we reduce to proving $h_n \not= 0$.

\meno {\em Step~6.\/}
Apply $E_{13}^n$ to~(\ref{sl3 step3}): since $E_{13} P_{13} = 0$, we obtain $E_{13}^n P_{23} (F_{13}^n) = E_{13}^n h_n F_{13}^n$.  Thus we reduce to proving $E_{13}^n P_{23} (F_{13}^n) \not= 0$.

\meno {\em Step~7.\/}
For any $H$ in $\dh$, define $d_k(H) := \prod_1^k (H + 1 + i)$.  By~(\ref{P sum}),
\begin{equation} \label{sl3 step7}
   E_{13}^n P_{23} (F_{13}^n) = 
   \sum_{k=0}^n\, \frac{(-1)^k / k!} {d_k(H_{23})^{-n \alpha_{13}}}\,
   E_{13}^n F_{23}^k E_{23}^k F_{13}^n,
\end{equation}
where $d_k(H_{23})^{-n \alpha_{13}}$ is the $-n \alpha_{13}$-shift of $d_k(H_{23})$ defined in~(\ref{shift}).  We are working in $M(\dsl_3)$, so~(\ref{sl3 step7}) is in $M(\dsl_3)_0$, which is $\dF(\dh)$.

We can conclude the proof efficiently with the following trick.  All of the denominators $d_k(H_{23})^{-n \alpha_{13}}$ with $k < n$ are strict divisors of $d_n(H_{23})^{-n \alpha_{13}}$, so if there is no cancellation between $d_n(H_{23})^{-n \alpha_{13}}$ and $E_{13}^n F_{23}^n E_{23}^n F_{13}^n$, then~(\ref{sl3 step7}) is non-zero.

Simplifying $E_{13}^n F_{23}^n E_{23}^n F_{13}^n$ in $M(\dsl_3)$, $E_{23}^n F_{13}^n$ becomes a $\bC$-multiple of $F_{12}^n$ and $E_{13}^n F_{23}^n$ becomes a $\bC$-multiple of $E_{12}^n$.  Hence the whole expression is a polynomial in $H_{12}$, which admits no cancellation with any polynomial in $H_{23}$.  $\Box$

\subsection{$\dsl_4$: Proof of Theorem~\ref{sl4}}\

\meno {\em Proof of Part~(i).\/}
By Lemma~\ref{PglPl}, the equation holds \iff\ multiplying its right side by $P_{12}$ gives $P_{1234}$.  Since $P_{12}$ is idempotent and commutes with the factors, we must prove that there is a unique $Q_{124}^{12}$ in $\dF(\dl_{124})^{\dl_{12}}$ satisfying
\begin{equation*}
   P_{1234} = (P_{123}^{12} P_{12}) Q_{124}^{12} (P_{34} P_{12})
   = P_{123} Q_{124}^{12} P_{12, 34}.
\end{equation*}

\meno {\em Step~1.\/}
By~(\ref{l decomp of M}), all $\dl_{123}$-HWVs in $M(\dsl_4)$ have weights in $-\Span_\bN \{\alpha_{14}, \alpha_{24}, \alpha_{34} \}$, and all $\dl_{12, 34}$-HWVs in $M(\dsl_4)$ have weights in $-\Span_\bN \{\alpha_{13}, \alpha_{14}, \alpha_{23}, \alpha_{24} \}$.  It follows that for any $Q_{124}^{12}$ in $\dF(\dl_{124})^{\dl_{12}}$, $P_{123} Q_{124}^{12} P_{12, 34}$ annihilates $M(\dsl_4)_{-\nu}$ for
\begin{equation*}
   \nu \not\in \Span_\bN \bigl\{\alpha_{14}, \alpha_{24}, \alpha_{34} \bigr\} \cap 
   \Span_\bN \bigl\{ \alpha_{13}, \alpha_{14}, \alpha_{23}, \alpha_{24} \bigr\} = \Span_\bN \bigl\{ \alpha_{14}, \alpha_{24} \bigr\}.
\end{equation*}

During this proof, write $n$ for an ordered pair $(n_{14}, n_{24})$ in $\bN^2$ and set
\begin{equation*}
   F^n := F_{14}^{n_{14}} F_{24}^{n_{24}}, \quad
   E^n := E_{14}^{n_{14}} E_{24}^{n_{24}}, \quad
   n \cdot \alpha := n_{14} \alpha_{14} + n_{24} \alpha_{24}, \quad
   |n| := n_{14} + n_{24}.
\end{equation*}
Equip $\bN^2$ with the usual partial order.  The preceding paragraph shows that it will suffice to choose $Q_{124}^{12}$ so that $P_{123} Q_{124}^{12} P_{12, 34}$ maps~$1$ to $1$ and annihilates $M(\dsl_4)_{-n \cdot \alpha}$ for $n > 0$.

\meno {\em Step~2.\/}
Using two PBW bases of $\dU(\dn^-)$, one with $F_{12}$ and $F_{34}$ to the left and the other with $F_{12}$, $F_{23}$, and $F_{13}$ to the left, we find that
\begin{equation} \label{sl4 step2}
   P_{123} \bigl( M(\dsl_4)_{-n \cdot \alpha} \bigr) = \dF(\dh) P_{123} (F^n), \quad
   P_{12, 34} \bigl( M(\dsl_4)_{-n \cdot \alpha} \bigr) = \dF(\dh) P_{12, 34} (F^n).
\end{equation}
By~(\ref{l decomp of M}), both of these spaces are non-zero.  By the second of the two equations, we are done if we prove that there is a unique choice of $Q_{124}^{12}$ such that $P_{123} Q_{124}^{12}$ maps $1$ to $1$ and $P_{12, 34}(F^n)$ to $0$ for $n > 0$.

\meno {\em Step~3.\/}
Note that $\Delta \bigl( \dU(\dl_{124}^-) \bigr) = - \Span_\bN \{ \alpha_{12}, \alpha_{24} \}$.  Lemma~\ref{nonstandard} implies that the sum of the weight spaces of $M(\dsl_4)$ with weights in this set is $\dl_{124}$-invariant and is a sum of copies of $M(\dl_{124})$, as follows:
\begin{equation*}
   \bigoplus_{\nu \in \Delta(\dU(\dl_{124}^+))} M(\dsl_4)_{-\nu} =
   \bigoplus_{j \in \bN^2} \dU(\dl_{124}^-) \dF(\dh) v_j, \quad
   v_{(j_{14}, j_{24})} := P_{124} (F_{13}^{j_{14}} F_{23}^{j_{24}} F_{34}^{|j|}).
\end{equation*}

Applying $P_{12}$ to this equation, we find that for $n \in \bN^2$ the space of $\dl_{12}$-HWVs in $M(\dsl_4)$ of weight $-n \cdot \alpha$ is
\begin{equation*}
   P_{12} \bigl( M(\dsl_4)_{-n \cdot \alpha} \bigr) =
   \Span_{\dF(\dh)} \bigl\{ P_{12} (F^i v_{n-i}):\, i \in \bN^2,\ 0 \le i \le n \bigr\}.
\end{equation*}
In particular, for unique elements $h_i$ of $\dF(\dh)$ with $0 \le i \le n$,
\begin{equation} \label{sl4 step3}
   P_{12, 34} (F^n) = \sum_{0 \le i \le n \in \bN^2} h_i P_{12} (F^i v_{n-i}).
\end{equation}

\meno {\em Step~4.\/}
By~(\ref{l decomp of M}), $M(\dl_{124})$ is the direct sum of the copies of $M(\dl_{12})$ with HWVs $P_{12}(F^k)$, where $k = (k_{14}, k_{24})$ runs over $\bN^2$.  Write $P_{124}^{12}[k]$ for the projection operator $P(\dl_{124}, \dl_{12}, F^k)$ projecting $M(\dl_{124})$ to $\dU(\dl_{12}) \dF(\dh) P_{12}(F^k)$ along the other summands of~(\ref{l decomp of M}).  Note that $P_{124}^{12}[k]$ is in $\dF(\dl_{124})^{\dl_{12}}$, and for any element $Q_{124}^{12}$ of $\dF(\dl_{124})^{\dl_{12}}$ there exist unique elements $q_k$ of $\dF(\dh)$ such that
\begin{equation} \label{sl4 step4}
   Q_{124}^{12} = \sum_{k \in \bN^2} q_k P_{124}^{12}[k].
\end{equation}

Now $P_{124}^{12}[k]$ acts on $P_{12} (F^i v_{n-i})$ by~$1$ if $i = k$ and by~$0$ otherwise.  Therefore
\begin{equation*}
   Q_{124}^{12} P_{12, 34} (F^n) = \sum_{0 \le i \le n \in \bN^2} q_i h_i P_{12} (F^i v_{n-i}).
\end{equation*}

\meno {\em Step~5.\/}
Apply $P_{123}$: because $P_{123} P_{12} = P_{123}$, we must choose the $q_k$ so that
\begin{equation} \label{sl4 step5}
   P_{123} Q_{124}^{12} P_{12, 34} (F^n) = \sum_{0 \le i \le n \in \bN^2 } q_i h_i P_{123} (F^i v_{n-i})
\end{equation}
is~$1$ for $n=0$ and~$0$ for all $n > 0$.  We choose them by induction on the partial order on $\bN^2$.  Clearly $q_0$ must be~$1$.  Suppose that $q_i$ has been determined for $i < n$.  By the first equation in~(\ref{sl4 step2}), every summand of~(\ref{sl4 step5}) is a multiple of $P_{123} (F^n)$.  Therefore there is a unique choice of $q_n$ such that the right side is zero \iff\ the coefficient of $q_n$ is non-zero.  Since $P_{123} (F^n) \not= 0$, we reduce to proving $h_n \not= 0$.

\meno {\em Step~6.\/}
Apply $E^n$ to~(\ref{sl4 step3}).  On the right side, the summand $E^n h_i P_{12} (F^i v_{n-i})$ is in the copy $\dU(\dl_{124}^-) \dF(\dh) v_{n-i}$ of $M(\dl_{124})$ generated by $v_{n-i}$.  But the weights of this space are all~$\le -(n - i) \cdot \alpha$, so only the summand at $i = n$ can be non-zero.  Thus $E^n P_{12, 34}(F^n)$ is $E^n h_n P_{12}(F^n)$, so we need only prove $E^n P_{12, 34}(F^n) \not= 0$.

\meno {\em Step~7.\/}
We use the same denominator trick used for $\dsl_3$.  Disregarding $\bC$-scalars, the summands of $P_{12, 34} = P_{12} P_{34}$ may be written as $F_{12}^a F_{34}^b E_{34}^b E_{12}^a$, with denominators $d_a(H_{12}) d_b(H_{34})$.  Consider
\begin{equation} \label{sl4 step7}
   E^n F_{12}^a F_{34}^b E_{34}^b E_{12}^a F^n =
   E_{24}^{n_{24}} E_{14}^{n_{14}} F_{12}^a F_{34}^b E_{34}^b E_{12}^a F_{14}^{n_{14}} F_{24}^{n_{24}}
\end{equation}
in $M(\dsl_4)_0 = \dF(\dh)$.  Because $E_{12}$ commutes with $F_{24}$, this term is zero for either $a > n_{14}$ or $b > |n|$.  Conversely, the denominators of the non-zero terms with either $a < n_{14}$ or $b < |n|$ are strict divisors of $d_{n_{14}}(H_{12}) d_{|n|}(H_{34})$.  Thus if~(\ref{sl4 step7}) is non-zero at $a = n_{14}$ and $b = |n|$ and admits no $\dF(\dh)$-cancellation against $d_{n_{14}}(H_{12}) d_{|n|}(H_{34})$, we are done.

Using $[E_{14}, F_{12}] = -E_{24}$ and $[E_{12}, F_{14}] = -F_{24}$, check that~(\ref{sl4 step7}) at $a = n_{14}$ and $b = |n|$ is $\bC$-proportional to $E_{24}^{|n|} F_{34}^{|n|} E_{34}^{|n|} F_{24}^{|n|}$.  This in turn is $\bC$-proportional to $E_{23}^{|n|} F_{23}^{|n|}$, which is a polynomial in $H_{23}$.  $\Box$

\meno {\em Proof of Part~(ii).\/}
By Lemma~\ref{PglPl}, the equation holds \iff\ multiplying the right side by $P_{23}$ gives $P_{1234}$.  Since $P_{23}$ is idempotent and commutes with the factors, we must prove that $P_{1234}$ is
\begin{equation*}
   (P_{123}^{23} P_{23}) Q_{14} (P_{234}^{23} P_{23})
   = P_{123} Q_{14} P_{234} = (P_{12} Q_{13} P_{23}) Q_{14} (P_{23} Q_{24} P_{34}).
\end{equation*}
This simplifies to $P_{12} Q_{13} P_{23} Q_{14} Q_{24} P_{34}$, which is $P_{1234}$ by Theorem~\ref{fin fac}.

We do not have a short proof of the uniqueness of $Q_{14}$.  An argument parallel to the one used to prove Part~(i) shows that there exists a unique element $\t Q_{14}$ of $\dF(\dl_{14})_0$ such that $P_{1234}^{23} = P_{123}^{23} \t Q_{14} P_{234}^{23}$; we will not give the details.  By the preceding paragraph, $\t Q_{14}$ must be $Q_{14}$.  $\Box$

\subsection{$\dsl_5$: Proofs of Theorems~\ref{sl5r1} and~\ref{sl5r2}}\

\meno {\em Proof of Theorem~\ref{sl5r1}~(i).\/}
Most of the arguments are similar to those in the proof of Theorem~\ref{sl4}~(i).  Multiplying by $P_{12}$ and applying Lemma~\ref{PglPl}, we come down to proving that there is a unique $Q_{125}^{12}$ in $\dF(\dl_{125})^{\dl_{12}}$ satisfying
\begin{equation} \label{sl5 step0}
   P_{12345} = P_{1234} Q_{125}^{12} P_{12, 345}.
\end{equation}

\meno {\em Step~1.\/}
By~(\ref{l decomp of M}), the weights of $P_{1234} \bigl( M(\dsl_5) \bigr)$ comprise $-\Span_\bN \{\alpha_{r, 5}: r < 5 \}$ and the weights of $P_{12, 345} \bigl( M(\dsl_5) \bigr)$ comprise $-\Span_\bN \{\alpha_{1, r}, \alpha_{2, r}: r > 2 \}$.  Conclude that for any $Q_{125}^{12}$ in $\dF(\dl_{125})^{\dl_{12}}$, $P_{1234} Q_{125}^{12} P_{12, 345}$ annihilates $M(\dsl_4)_{-\nu}$ for
\begin{equation*}
   \nu \not\in \Span_\bN \bigl\{ \alpha_{15}, \alpha_{25} \bigr\}.
\end{equation*}

During this proof, write $n$ for an ordered pair $(n_{15}, n_{25})$ in $\bN^2$ and set
\begin{equation*}
   F^n := F_{15}^{n_{15}} F_{25}^{n_{25}}, \quad
   E^n := E_{15}^{n_{15}} E_{25}^{n_{25}}, \quad
   n \cdot \alpha := n_{15} \alpha_{15} + n_{25} \alpha_{25}, \quad
   |n| := n_{15} + n_{25}.
\end{equation*}
It will suffice to choose $Q_{125}^{12}$ so that $P_{1234} Q_{125}^{12} P_{12, 345} \bigl( M(\dsl_5)_{-n \cdot \alpha} \bigr) = \delta_{0,n} \dF(\dh)$.

\meno {\em Step~2.\/}
Using appropriate PBW bases of $\dU(\dn^-)$, deduce that
\begin{equation*}
   P_{1234} \bigl( M(\dsl_5)_{-n \cdot \alpha} \bigr) = \dF(\dh) P_{1234} (F^n), \quad
   P_{12, 345} \bigl( M(\dsl_5)_{-n \cdot \alpha} \bigr) = \dF(\dh) P_{12, 345} (F^n).
\end{equation*}
By~(\ref{l decomp of M}), both spaces are non-zero.  We are done if we prove that there is a unique choice of $Q_{125}^{12}$ such that $P_{1234} Q_{125}^{12}$ maps $P_{12, 345}(F^n)$ to $\delta_{0,n}$.

\meno {\em Step~3.\/}
By Lemma~\ref{nonstandard}, $\bigoplus_{\nu \in \Delta(\dU(\dl_{125}^+))} M(\dsl_5)_{-\nu}$ is $\dl_{125}$-invariant and is a sum of copies of $M(\dl_{125})$.  As in the case of $\dsl_4$, the highest weights of these copies of $M(\dl_{125})$ are all in $-\Span_\bN \{ \alpha_{15}, \alpha_{25} \}$, but here there is more than one copy for each highest weight.  Writing temporarily $\I_{125}$ for the set of $I$ in $\I(\dsl_5, \dl_{125})$ such that $|I|$ is in $\Span_\bN \{ \alpha_{15}, \alpha_{25} \}$, Lemma~\ref{nonstandard} leads to
\begin{equation*}
   \bigoplus_{\nu \in \Delta(\dU(\dl_{125}^+))} M(\dsl_5)_{-\nu} =
   \bigoplus_{I \in \I_{125}} \dU(\dl_{125}^-) \dF(\dh) P_{125} (F^I).
\end{equation*}

Because $P_{12} \bigl( M(\dl_{125}) \bigr)$ has 1-dimensional weight spaces, we find that for $0 \le j \le n$ in $\bN^2$ there are unique $\dl_{125}$-HWVs $v_j$ of weight $-j \cdot \alpha$ in $M(\dsl_5)$ such that
\begin{equation} \label{sl5 step3}
   P_{12, 345} (F^n) = \sum_{0 \le i \le n \in \bN^2} P_{12} (F^i v_{n-i}).
\end{equation}
In this notation, the coefficients $h_i$ in~(\ref{sl4 step3}) have been absorbed by the $v_{n-i}$.  The key point is that $v_0$ is in $\dF(\dh)$; we will not be concerned with the other $v_j$.

\meno {\em Step~4.\/}
The first paragraph remains as for $\dsl_4$ except that $4$ is replaced by~$5$:
\begin{equation*}
   Q_{125}^{12} = \sum_{k \in \bN^2} q_k P_{125}^{12}[k].
\end{equation*}
Since $P_{125}^{12}[k] \bigl( P_{12} (F^i v_{n-i}) \bigr) = \delta_{i, k} P_{12} (F^i v_{n-i})$, we have
\begin{equation*}
   Q_{125}^{12} P_{12, 345} (F^n) = \sum_{0 \le i \le n \in \bN^2} q_i P_{12} (F^i v_{n-i}).
\end{equation*}

\meno {\em Step~5.\/}
Apply $P_{1234}$: because $P_{1234} P_{12} = P_{1234}$, we must choose the $q_k$ so that
\begin{equation*}
   P_{1234} Q_{125}^{12} P_{12, 345} (F^n) = \sum_{0 \le i \le n \in \bN^2 } q_i P_{1234} (F^i v_{n-i})
   = \delta_{0, n}.
\end{equation*}
Proceed as for $\dsl_4$: there is a unique solution for $q_n$ \iff\ $v_0 \not= 0$.

\meno {\em Step~6.\/}
Apply $E^n$ to~(\ref{sl5 step3}): by the weight argument used before, the right side becomes $E^n P_{12}(F^n) v_0$, so we need only prove $E^n P_{12, 345}(F^n) \not= 0$.

\meno {\em Step~7.\/}
Recall that $P_{12, 345} = P_{12} P_{345}$ and $P_{345} = P_{45} P_{35} P_{34}$.  Since $P_{34} (F^n) = F^n$, the denominators of the summands of $E^n P_{12, 345}(F^n)$ are polynomials in $H_{12}$, $H_{35}$, and $H_{45}$, and the numerators are terms like
\begin{equation*}
   E_{25}^{n_{25}} E_{15}^{n_{15}} F_{12}^a F_{45}^b E_{45}^b
   F_{35}^c E_{35}^c E_{12}^a F_{15}^{n_{15}} F_{25}^{n_{25}}.
\end{equation*}
The largest non-zero numerator occurs at $a = n_{15}$ and $b = c = |n|$ and simplifies in $M(\dsl_5)$ to a polynomial in $H_{23}$, so the result follows as before.  $\Box$

\meno {\em Proof of Theorem~\ref{sl5r1}~(ii).\/}
Arguing as for $\dsl_4$, the equation holds \iff\ multiplying the right side by $P_{23}$ gives $P_{12345}$, and so we come down to proving that $P_{12345}$ is $P_{1234} Q_{15} P_{2345}$.  Applying Theorem~\ref{fin fac} to $P_{12345}$, $P_{1234}$, and $P_{2345}$ confirms this.  Again, we do not have a short proof of uniqueness: it is necessary to follow the steps leading to Part~(i).  $\Box$

\meno {\em Proof of Theorem~\ref{sl5r2}.\/}
We will omit the proof of Part~(i): it is similar to the proofs of Theorem~\ref{sl4}~(i) and Theorem~\ref{sl5r1}~(i), using $F^n := F_{15}^{n_{15}} F_{25}^{n_{25}} F_{35}^{n_{35}}$ in place of the earlier definitions of $F^n$.  Regarding Part~(ii), we mention only that multiplying by $P_{234}$ and following the argument for Theorem~\ref{sl5r1}~(ii) proves the formula; for uniqueness we must go through all seven steps.  On the other hand, multiplying by $P_{12, 34}$ in Part~(iii) brings us to the proof of Theorem~\ref{sl5r1}~(i) at~(\ref{sl5 step0}), proving the result completely with no additional work.  $\Box$

\section{Remarks} \label{Remarks}

\subsection{}
Our primary goal at this point is to characterize those $\dg$ and $\dl$ for which our method gives a non-commutative factorization of $P(\dg, \dl)$.  It does not apply in multiply laced cases, even for $\dl = \dh$: it is instructive to examine its failure for $P(\dgo_5)$.  It can also break down when $\dl_\ss$ is not simple; this type of failure first occurs for $\dg = \dsl_5$ and $\dl = \dl_{12, 45}$, the ``missing case'' in Theorem~\ref{sl5r2}.  Here the natural conjecture is that there is a unique element $Q_{1245}^{12, 45}$ of\/ $\dF(\dl_{1245})^{\dl_{12, 45}}$ such that
\begin{equation*}
   P_{12345}^{12, 45} = P_{123}^{12}\, Q_{1245}^{12, 45}\, P_{345}^{45}.
\end{equation*}
Multiplying by $P_{12, 45}$ and applying Lemma~\ref{PglPl}, we may replace this equation by
\begin{equation*}
   P_{12345} = P_{123, 45}\, Q_{1245}^{12, 45}\, P_{12, 345}.
\end{equation*}

Carrying out Step~1, we find that we need only consider $M(\dsl_5)_{-\nu}$ for
\begin{equation*}
   \nu \in \Span_\bN \bigl\{ \alpha_{14}, \alpha_{24}, \alpha_{15}, \alpha_{25} \bigr\}.
\end{equation*}
However, when we apply $P_{123, 45}$ and $P_{12, 345}$ to these $M(\dsl_5)_{-\nu}$ in Step~2, we do not in general obtain 1-dimensional images over $\dF(\dh)$.  For example, $F_{14} F_{25}$ and $F_{15} F_{24}$ have the same weight but, by Lemma~\ref{nonstandard}, independent images.  This causes Step~5 to fail.  Maybe the failure can be repaired using the fact that here the projections $P_{1245}^{12, 45}[k]$ appearing in Step~4 have more than one $k$ corresponding to each weight, but we have not yet overcome the difficulty.

\subsection{}
We would also like to have some description of the factors $Q(\dm, \dl)$ of $P(\dg, \dl)$.  An explicit formula would be best, but short of that one could try to prove that they have certain properties possessed by the factors $Q_{\rho_\dg(H_\alpha)}(\da_\alpha)$ of $P(\dg)$.

For example, consider~(\ref{sl3 step4}).  As mentioned below~(\ref{Pt fac}), in fact only $q_0$ and $q_1$ are non-zero, and in the analogous expression for $Q_t(\dsl_2)$, only $q_0, \ldots, q_{t-1}$ are non-zero.  It is natural to predict that this phenomenon occurs also in the relative case.  The first instance of this is~(\ref{sl4 step4}), where we expect that only the first few $q_k$ are non-zero: probably only $q_{(0,0)}$, $q_{(1,0)}$, and $q_{(0,1)}$.  Relative versions of~(\ref{Pt fac}) giving infinite commutative factorizations of the $Q(\dm, \dl)$ would resolve the situation.

To give another example, recall that for $\dl = \dh$, the nonstandard subalgebras $\dm$ corresponding to the AST factors of $P(\dg)$ are $\dl_\alpha = \da_\alpha + \dh$.  As pointed out in Section~\ref{sl3 proof}, in this setting our method only proves the existence of factors $Q(\dl_\alpha, \dh)$ in $\dF(\dl_\alpha)^\dh$, although in fact they are in $\dF(\da_\alpha)^\dh$.  Does this generalizes to all $\dl$?  That is, are the factors $Q(\dm, \dl)$ in $\dF(\dm_\ss)^\dl$?  In the context of the preceding paragraph, affirming this amounts to proving that the coefficients $q_k$ are in $\dF(\dh \cap \dm_\ss)$.

\subsection{}
We conclude by explaining some implications of our results for Conjecture~\ref{rel conj}.  In all the cases we treat, the subalgebras $\dm$ and the factors $Q(\dm, \dl)$ have the following properties:

\begin{enumerate}

\smallbreak \item[(i)]
$\Delta(\dn^+) \backslash \Delta(\dl^+)$ is the disjoint union $\coprod_\dm \Delta(\dm^+) \backslash \Delta(\dl^+)$.

\smallbreak \item[(ii)]
$Q(\dm, \dl)$ is the projector $P(\dm, \dl)$ \iff\ $\dm$ is standard in $\dg$.

\end{enumerate}
\smallbreak
We expect that these properties hold in general.  When they do, Conjecture~\ref{rel conj} would follow immediately if one could prove that there is a formula for $Q(\dm, \dl)$ with total denominator
\begin{equation} \label{Q denom}
   D(\dm, \dl) := \prod_{i=1}^\infty \, \HC_\dl^{-1} \,
   \biggl( \, \prod_{\alpha \in \Delta(\dm^+) \backslash \Delta(\dl^+)} \,
   \bigl(H_\alpha + i \bigr)^{\rho_\dg} \biggr).
\end{equation}

In several cases in Theorems~\ref{sl4}, \ref{sl5r1}, and~\ref{sl5r2}, $\dm = \dl + \da_\alpha$ for some $\alpha$ in $\Delta(\dn^+)$ such that $\da_\alpha$ commutes with $\dl$.  In all of these cases the factor $Q(\dm, \dl)$ is simply $Q_{\rho_\dg (H_\alpha)} (\da_\alpha)$ from~(\ref{Pt sum}).  For such $\dm$, $H_\alpha$ is in $\dz(\dl)$, so~(\ref{Q denom}) does indeed reduce to the denominator of $Q(\dm, \dl)$.  Coupling these observations with Theorem~\ref{13} proves Conjecture~\ref{rel conj} for $P(\dsl_4, \dl_{23})$ and $P(\dsl_5, \dl_{234})$.  More generally, using the proof of Theorem~\ref{sl4}~(ii) one easily obtains:

\begin{lemma}
For $1 < a < b < n$, $P_{1 \cdots n}^{a\, a+1 \cdots b} = P_{1 \cdots n-1}^{a\, a+1 \cdots b} Q_{1 n} P_{2 \cdots n}^{a\, a+1 \cdots b}$.
\end{lemma}

Hence Theorem~\ref{13} and induction on $n$ give:

\begin{prop}
For $1 \le a < b \le n$, Conjecture~\ref{rel conj} holds for $P(\dsl_n, \dl_{a\, a+1 \cdots b})$.  Put differently, the conjecture holds if\/ $\dg$ is of type $A_n$ and the simple roots of\/ $\dl$ form a connected segment of the Dynkin diagram of\/ $\dg$.

\end{prop}

\bibliographystyle{amsalpha}

\end{document}